\documentclass[a4paper,10pt]{article} 
\usepackage[english]{babel} 
\usepackage[matrix,arrow,tips,curve]{xy}
\usepackage{amsmath}
\usepackage{amssymb,amsthm}
\usepackage{enumerate}

\usepackage[dvips]{graphics}

\usepackage{psfrag}

\newtheorem{thm}{Theorem}
\newtheorem{lemma}[thm]{Lemma}
\newtheorem{prop}[thm]{Proposition}
\newtheorem{cor}[thm]{Corollary}
\newtheorem*{fact}{Fact}

\theoremstyle{definition}
\newtheorem{defi}[thm]{Definition}

\theoremstyle{remark}
\newtheorem*{remark}{Remark}
\newtheorem*{claim}{Claim}
\newtheorem{example}{Example}

\newcommand{\g}{\Gamma}

\newcommand{\Z}{\mathbb{Z}}

\newcommand{\C}{\mathbb{C}}

\newcommand{\Pic}{\operatorname{Pic}}

\newcounter{prg}[section]

\newcommand{\la}{\longrightarrow}
\def\ev{\Delta}
\def\Mg{\overline{M_g}}
\def\Sg{\overline{S_g}}
\def\Pdg{\overline{P_{d,\,g}}}
\def\Pgg{\overline{P_{g-1,\,g}}}
\def\Yn{Y^{\nu}}
\def\Ln{L^{\nu}}
\def\X{\mathcal X}
\def\O{\mathcal O}

\psfrag{uno}{\fontsize{20}{20}$B_{\g}=\{0,1\}$}
\psfrag{due}{\fontsize{20}{20}$B_{\g}=\{0,1,2\}$}
\psfrag{W1}{\fontsize{20}{20}$W_1$}
\psfrag{W2}{\fontsize{20}{20}$W_2$}
\psfrag{Case1}{\fontsize{100}{100}{\huge Case (1)}}
\psfrag{Case2}{\fontsize{60}{60}{\huge Case (2)}}
\psfrag{C1}{\fontsize{20}{20}$C_1$}
\psfrag{C2}{\fontsize{20}{20}$C_2$}
\psfrag{C3}{\fontsize{20}{20}$C_3$}
\psfrag{C}{\fontsize{20}{20}$C$}

\begin{document}
\title{Combinatorial properties of stable spin curves}
\author{Lucia Caporaso and Cinzia Casagrande}
\date{ }
\maketitle

\bigskip

\noindent Lucia Caporaso \\
Dipartimento di Matematica, Universit\`a di Roma Tre \\
Largo S.\ L.\ Murialdo, 1 \\
00146 Roma - ITALY\\
caporaso@mat.uniroma3.it

\bigskip

\noindent Cinzia Casagrande \\
Dipartimento di Matematica, Universit\`a di Roma ``La Sapienza" \\
Piazzale Aldo Moro, 2 \\
00185 Roma - ITALY\\
ccasagra@mat.uniroma1.it

\bigskip

\

\hfill Dedicated to Professor Steven Kleiman, for his birthday


\section{Stable spin curves and their moduli}
 \subsection {Summary}
\nocite{kleiman2}
The first part of this paper\footnote{2000 Mathematics Subject
Classification: 14H10, 05C75} 
describes the  moduli space of stable spin curves (constructed in \cite{cornalba1}),
explaining how its geometry   is governed by the
combinatorics of stable curves.  In this context, the standard graph theoretic framework
(where to every stable curve one associates its ``dual graph") is not just a book-keeping device:
in Section 2, some purely combinatorial results are proved
 (Theorems \ref{due} 
and \ref{tre}), having moduli theoretic applications.
More precisely, certain strata of the moduli space of stable curves  are characterized by a (finite)
set of integers that  measures the non-reducedness of the  scheme of spin curves,
and which is definable in purely graph-theoretical terms (Definition \ref{betti}).

 \subsection{The basic functors and their compactifications}
A smooth spin curve is a pair $(X,L)$ where $X$ is a smooth, connected, projective curve of genus $g$
and $L$ a
theta-characteristic of $X$, that is, a line bundle such that $L^{\otimes 2}=K_X$.
Smooth spin curves naturally define a coarsely representable functor: 
to any  family 
$\X \la B$  of smooth curves (i.e. the fiber $X_b$ over every point $b\in B$ is a smooth curve) it
associates the relative spin curve 
$$\mathcal{S} _{\mathcal X/B}\la B$$ 
whose fiber over $b\in B$ is the set of
$2^{2g}$  smooth spin curves supported on $X_b$
(i.e. the set of theta-characteristics of $X_b$).
To complete the picture, consider the Picard functor: denote by 
$$
\Pic ^d_{{\mathcal {X}}/B}\la B
$$
the relative, degree $d$, Picard variety, whose fiber over $b\in B$ is the variety parametrizing line
bundles of degree $d$ on $X_b$.
We have a commutative diagram:
$$
\xymatrix{{{\mathcal{S}} _{{\mathcal
        {X}}/B}\,}\ar[dr]\ar@{^{(}->}[r]&{\Pic ^{\,g-1}_{{\mathcal
        {X}}/B}}\ar[d]\\&{B}}
$$
Because of the coarse representability of all the functors involved, we get a global picture
\begin{equation}
\label{diagram}
\xymatrix{{S_g\,}\ar[dr]\ar@{^{(}->}[r]&{P_{g-1,g}}\ar[d]\\&{M_g}}
\end{equation}
where $M_g$ is the moduli space of smooth curves of genus $g\geq 2$, $S_g$ the moduli space of smooth spin curves and 
$P_{d,g}$ the universal Picard variety of degree $d$, parametrizing isomorphism classes of 
degree $d$ line bundles over smooth curves of genus $g$ (the diagram, of course,  represents
only the case
$d=g-1$).

The word ``space" here means either algebraic stack or algebraic scheme; the warning is the usual one: in
the category of schemes, the moduli properties of the above diagram  fail for objects with non-trivial
automorphisms.

In order to compactify $M_g$, P. Deligne and D. Mumford introduced {\it stable curves}  (in 
\cite{dm}). 

\begin{defi}
\label{stable}
A stable curve  is a reduced, connected curve having only
ordinary double points as singularities and ample dualizing sheaf.
This last condition is equivalent to the fact that every smooth rational component of the curve
contains at least $3$ nodes.

By weakening this last requirement, asking that the nodes contained in any smooth rational component be at least $2$,
one gets the definition of a {\it semistable curve}. Thus a semistable curve $X$ is a nodal curve that fails
from being stable if it contains some smooth, rational component $E$  such that 
$\# ( {\overline{X\smallsetminus E}}\cap X)=2$. Such an $E$ will be
called {\it exceptional}. 

A {\it quasistable curve} is a semistable curve such that two exceptional components never meet each
other.

The {\it stable model}
of a  semistable  curve is the uniquely defined stable curve obtained by contracting every exceptional
component to a point.
\end{defi}
Stable curves of (arithmetic) genus $g\geq 2$ have a moduli space, denoted by $\Mg$, which is projective and contains
$M_g$ as a dense open subspace.

 Having compactified $M_g$ in this fashion, the problem arises on how to compactify $S_g$ and
$P_{d,g}$ accordingly.

A solution for $S_g$ was  given by M. Cornalba in 
\cite{cornalba1}, consistently with the Deligne-Mumford
construction. He defined  {\it stable spin  curves}:
\begin{defi}
\label{spin} A stable spin curve is a pair
$(Y,L)$ where $Y$ is a quasistable curve and $L$ a line bundle on $Y$ with the following properties.
Denote by $E$ any exceptional component of $Y$ and by $Z:=\overline{Y\smallsetminus\cup E}$ the closure of the
complement of
all exceptional components; then  the restriction  of $L$ to every exceptional component $E$  is
${\O }_E(1)$,  and the restriction to what remains satisfies
$$
L_{|Z}^{\otimes 2}\cong \omega _Z.
$$
\end{defi}
Notice that the 
degree of $L$ is $g-1$.

Stable spin curves are shown in \cite{cornalba1} to  have a projective moduli space $\Sg$, with a natural,
finite morphism of degree $2^{2g}$ onto $\Mg$
$$
\pi\colon\Sg \la \Mg .
$$
The fiber of $\pi$ over a stable curve $X\in \Mg$ is a zero-dimensional scheme
parametrizing  stable spin curves $(Y,L)$ such that the stable model of $Y$ is $X$.

With diagram \eqref{diagram} in mind, the question remains on how to compactify $P_{d,g}$; 
the problem of completing the Picard functor is rich with many  aspects,
which will not be described here, and has 
been (and still is) the object of interest for
a long time. 
For  solutions and methods to approach it, that differ from what will be presented in this paper, 
we refer to the recent articles \cite{ak} and \cite{ekg} of  A. Altman, S.
Kleiman, E. Esteves, M. Gagn\'e, and to the references
therein.

In this paper, we  are interested in stable curves
(thus, we only allow  nodal singularities); 
furthermore, we consider completions over
$\Mg$ of the Spin functor
and, marginally, of the Picard functor, that use polarized quasistable curves as boundary points.
The first to be completed  was $S_g$, by
the above described  space $\Sg$.

A   compactification, $\Pdg$,  of $P_{d,g}$ over $\Mg$,
 was later constructed in \cite{caporaso}; briefly said, $\Pdg$ parametrizes pairs $(Y,M)$ where $Y$ is a quasistable
curve of genus $g$, $M$ is a line
bundle of degree $d$ on $Y$, having degree $1$ on all exceptional components
and satisfying other ``degree constraints" (which we shall not explain here).    

The evident similarity between the boundary points of $\Sg$ and $\Pdg$ seemed a bit striking,
since the  two constructions were  independent and used  different techniques.
Only recently,  C.~Fontanari (\cite{fontanari}) showed this
analogy is not an accident: 
he proves that $\Sg$ is naturally a subscheme of $\Pgg$, so that one simultaneously
compactifies  all objects in  \eqref{diagram}:
 $$\xymatrix{{\Sg\,}\ar[dr]\ar@{^{(}->}[r]&{\Pgg}\ar[d]\\&{\Mg}}$$
The above diagram, where the standard functoriality properties are satisfied, clarifies and highlights the ``naturality"
of the boundary objects: stable curves, spin stable curves, suitably polarized quasistable curves.
We shall keep it in mind throughout the  paper.

 \subsection{Stable spin curves}
We start by  recalling some  facts about stable spin curves, referring to
\cite{cornalba1} for detailed proofs.
Fix a stable curve  $X$ and  let $(Y,L)$ be a stable spin curve such that the stable model of $Y$ is $X$; the quasistable curve $Y$ will
be called the {\it {support}} of the spin curve.
We shall denote by $S_X$ the moduli space of stable spin curves whose support has $X$  as
stable model          . Thus, if $X$ has trivial automorphism group,  $S_X$ is the
the (scheme-theoretic) fiber of $\pi:\Sg\la \Mg$ over  $X$.

Let $\nu\colon 
X^{\nu} \la X$ be the normalization map, $C\subset X$ an irreducible component
and $C^{\nu}$ the corresponding component in $X^{\nu}$.
Let $\ev\subset X_{sing}$ be a 
set of nodes of $X$, set $\ev _C :=\ev \cap C$ and denote by
$D_C\subset C^{\nu}$ the preimage of $\ev _C$, i.e.
$D_C:=\nu ^{-1}(\ev _C)$.
Thus $D_C$ is an effective, reduced  divisor of $C^{\nu}$.
\begin{defi}
\label{even}
We say that $\ev$ is {\it even} if, for every irreducible component $C$ of $X$, $\deg D_C$ is even.
\end{defi}
For example, the empty set is even.

Notice that, equivalently, $\ev$ is even if, denoting by $Z$ the partial normalization of $X$ at all nodes that
are not in $\ev$, the dualizing sheaf $\omega _Z$ has even degree on every irreducible component of $Z$.

Consider now the set of all quasistable curves having $X$ as stable model;
  this set is obviously finite and in bijective correspondence with the set of subsets of
nodes of $X$. More precisely,  
let $Y$ be a quasistable curve 
and  let $\sigma\colon Y\la X$ be the natural morphism contracting all the exceptional components of $Y$.
 Denote by $\ev _Y\subset X_{sing} $ the set of nodes corresponding to the nodes of $Y$ that are not contained in
an exceptional
component:
$$
\ev _Y := \sigma (Y\smallsetminus\cup E)_{sing}.
$$
Clearly  $Y$  uniquely determines $\ev _Y$ and, conversely, for every $\ev \subset X_{sing}$ there exists a unique
quasistable curve $Y$ such that $\ev _Y = \ev$.

 A basic consequence of Cornalba's construction is the following
\begin{fact}
A quasistable curve  $Y$ is the support of a spin curve  if and only if $\ev _Y$ is even.
\end{fact}

Having characterized all quasistable curves appearing as supports of
spin curves, we fix one, $Y$, and 
describe all line bundles $L\in \Pic Y$ such that $(Y,L)\in S_X$; as always,
$\nu\colon\Yn \la Y$ denotes the normalization. 

Denote $L^{\nu} := \nu ^* L$, so that $L^{\nu}$ is the datum of a line bundle on each
irreducible component of
$\Yn$. By Definition \ref{spin}, we have 
$$
\Ln _{|E} = {\mathcal O}_E(1)
$$
and, for every non-exceptional component $C$ of $Y$
\begin{equation}
\label{C}
(\Ln _{|C^{\nu}})^{\otimes 2} = K_{C^\nu}\otimes{\mathcal O}(D_C).
\end{equation}
This last formula follows from Definition \ref{spin}: recall that, 
denoting $Z:=\overline{Y\smallsetminus\cup E}$, $L$ satisfies
\begin{equation}
\label{ZZ}
L_{|Z}^{\otimes 2} \cong \omega _Z 
\end{equation}
Fix now $\Ln$ as above; the 
set of all  line bundles $L$  on $Z$ that satisfy (\ref{ZZ}), and pull back to $\Ln$ (restricted to $Z^{\nu}$), is found
by looking at 
  the  exact sequence of algebraic groups:
\begin{equation}
\label{Z}
1\la (\C ^*)^{b_1(\Gamma _Z)}\la \Pic Z \xrightarrow{\ \nu^*} \Pic Z^{\nu} \la 0
\end{equation}
where $\Gamma _Z$ is the dual graph of $Z$ (whose definition we recall below), and $b_1(\Gamma
_Z)$ its first Betti number.
\begin{defi} Let $Z$ be a reduced nodal curve. 
The {\it dual graph} of $Z$, denoted by $\Gamma _Z$,
is the graph whose vertices are the irreducible components of $Z$ and whose edges are the nodes of $Z$.
\end{defi}
The above sequence says that there are $(\C ^*)^{b_1(\Gamma _Z)}$ line bundles on $Z$, all pulling back to
$\Ln_{|{Z^{\nu}}}$.  Of these, there are exactly $2^{b_1(\Gamma _Z)}$ line bundles 
that satisfy (\ref{ZZ}), in fact, on every node of $Z$, there are exactly two   gluings compatible with (\ref{ZZ}).

The gluing data on the remaining nodes of $Y$, lying on some exceptional component,  do not give different isomorphism classes of
spin
curves.  More precisely, let $N$ be any node of $Y$ that lies on an
exceptional component; different gluings of $\Ln$ over $N$ give the same  point of $S_X$,
but determine the scheme structure of $S_X$ at such a point, being responsible for the non-reducedness
of $S_X$.

Set
$$b=b_1(\Gamma _X) =b_1(\Gamma _Y)$$ 
and notice that
$$
 b=  \# Y_{sing} - \#\{\text{irreducible components of } Y\} +1 =g-p
$$
where 
$$p:= \sum _Cp_g(C),$$ 
$p_g(C)$ is the geometric genus of $C$ and the
the sum  is extended to all irreducible components $C$ of $X$ (or of $Y$). 

By what we said,  $S_X$ contains $2^{b_1(\Gamma _Z)}$ distinct points $(Y,L)$ such that $\nu ^*L=\Ln$.
At each of these points, the multiplicity of $S_X$ is computed in
\cite{cornalba1}, Section 4 and in \cite{cs}
Section 2.2; it is equal to $2^{b-b_1(\Gamma _Z)}$.

Let us recapitulate and check the above analysis by computing the length of $S_X$ (which must, of course, be equal to
$2^{2g}$).

The number of choices for $Y$ is equal to the number of even subsets on nodes of $X$, which is equal to the number of
 cyclic subgraphs of $\Gamma _X$ (see  the remark  in Section \ref{cycles}), which is equal to $2^b$.

The number of choices of $\Ln \in \Pic \Yn$ is $2^{2p}$ (recall that $p=\sum p_g(C)$);
in fact, for every irreducible component $C$ of $Y$, the number of
choices for the restriction of $\Ln$ to $C^{\nu}$ is equal to
$2^{2p_g(C)}$ (by formula (\ref{C})). 

For each $\Ln$ we have $2^{b_1(\Gamma _Z)}$ distinct points of $S_X$, 
corresponding to the different gluings over the nodes of $Z$ (by
(\ref{Z})). All such points  have multiplicity $2^{b-b_1(\Gamma _Z)}$. 
Summarizing:
$$
\text{length}(S_X)=2^b\cdot 2^{2p}\cdot 2^{b_1(\Gamma _Z)}\cdot 2^{b-b_1(\Gamma _Z)} = 2^{2b+2p}=2^{2g}.
$$

We shall now make a small change of notation, to better highlight how
the   scheme structure of $S_X$  is  governed by the combinatorics
of $X$.

To a set $\ev$ of nodes of $X$, one   associates  a subgraph of $\Gamma _X$ as follows:
 $\ev$  is identified with  a set of edges of $\Gamma _X$, which naturally generates a subgraph,
which is the smallest subgraph of $\Gamma _X$ containing $\ev$. For example, in the above set up,
the graph associated to $\ev _Y$ is identifiable to the dual graph of $Z$, $\Gamma _Z$.
We will abuse notation and denote by the same symbol,
$\ev$, the set of nodes of the curve, the corresponding set of edges of the dual graph,
and the subgraph generated by such edges.
Therefore, in the previous set up
$$
b_1(\Gamma _Z)= b_1(\ev _Y).
$$
Having done that,
we summarize the above analysis of the  structure of $S_X$:
\begin{prop} [Numerics of $S_X$]
\label{num} $S_X$ is a zero-dimensional scheme of length $2^{2g}$. The number of its
irreducible components is
$$
2^{2p}\cdot \Bigl(\sum _{\ev \subset X_{sing},\, \ev \ even}2^{b_1(\ev )}\Bigr)
$$
A component of $S_X$ parametrizing the stable spin curve $(Y,L)$ appears with multiplicity equal to
$2^{b-b_1(\ev _Y)}$.
\end{prop}

\subsection{Examples}
\begin{example}[$S_X$  is reduced if and only if $X$ is of compact type]
\label{ct}
Let $X$ be a curve of compact type;
then $b=0$ and the only even set of nodes of $X$ is $\emptyset$. Obviously, $b_1(\emptyset)=0$,
therefore $S_X$ is reduced and all of its points parametrize pairs $(Y,L)$ such that $Y$ is the quasistable curve
obtained by ``blowing up'' 
every node of $X$ and $L$ is the datum of a theta characteristic on
every irreducible component of $X$.

Conversely, suppose that $S_X$ is reduced. Then, by the proposition, for every even subset $\ev$ of nodes of $Y$ we must
have
$$b-b_1(\ev )=0.
$$
In particular ($\emptyset $ is even) 
$$b-b_1(\emptyset )=b=0
$$
therefore $X$ is of compact type.
\end{example}

\begin{example}[b=1]
\label{b1}
If $b=1$ then $p=g-1$. Denote by $W$  a (possibly empty) disjoint union of  a finite number of curves of compact type
$W_i$:
$W=\cup W_i$, with $W_i\cap W_j =\emptyset$.
There are two possibilities: 
\begin{enumerate}[(a)]
\item
  $X=C\cup W$ with $W$   as above, $C$  irreducible with one node $N$,
and   $\# C\cap W_i =1$ for every $i$;
\begin{center}  

 \scalebox{0.4}{\includegraphics{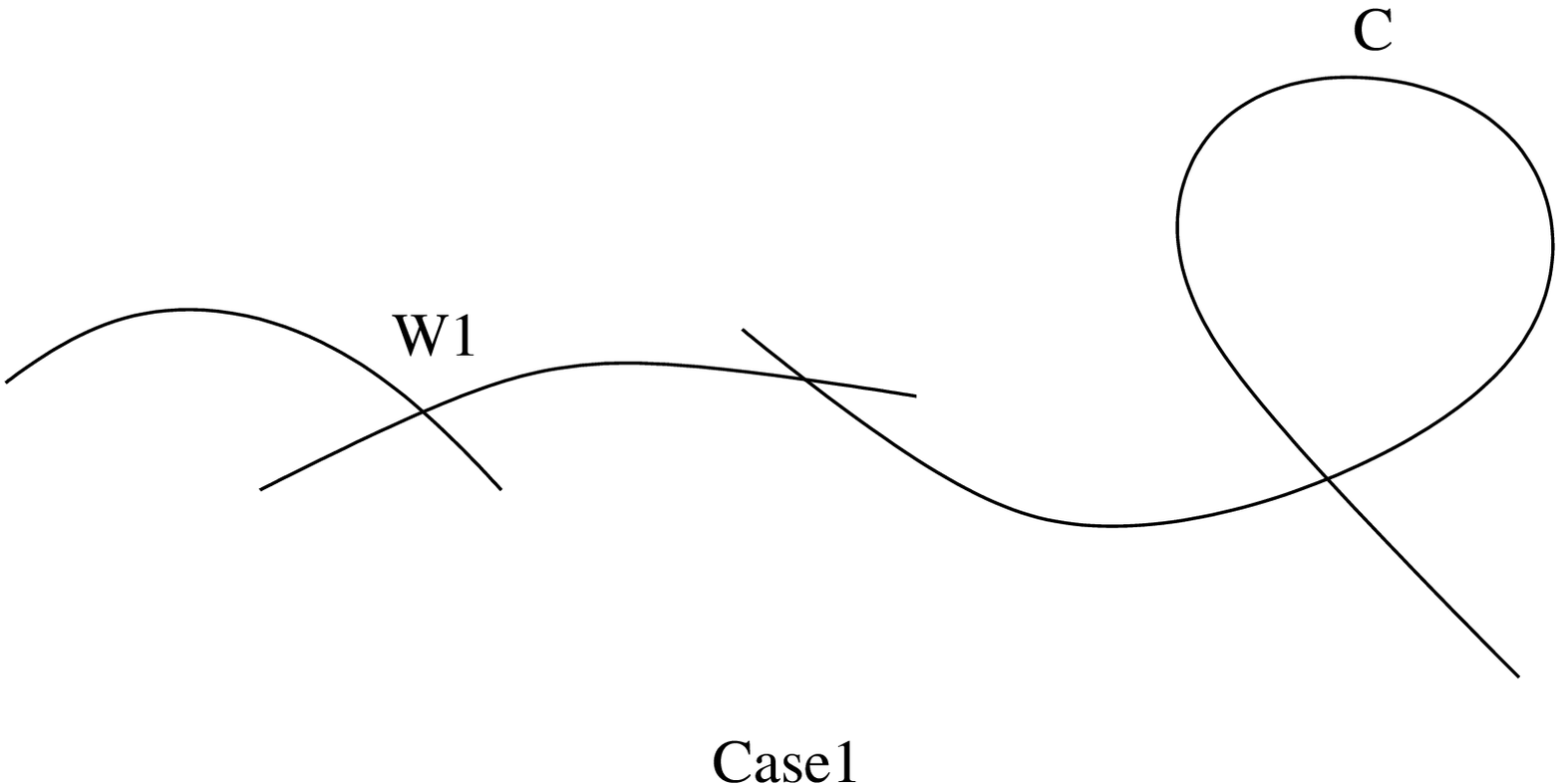}}

\end{center}
\item
 $X=D\cup W$ such that  $W$  is  as above, 
$D$ is a 
``cycle" of smooth components, that is, 
$D=C_1\cup\ldots \cup C_h $
with $C_i$ smooth, $C_i$ intersects $C_j$ if and only if $i$ and $j$ are
consecutive integers, or $i=1$ and $j=h$; finally 
 $\# (W_i\cap D)=1$ for every $i$.
\begin{center}  

 \scalebox{0.4}{\includegraphics{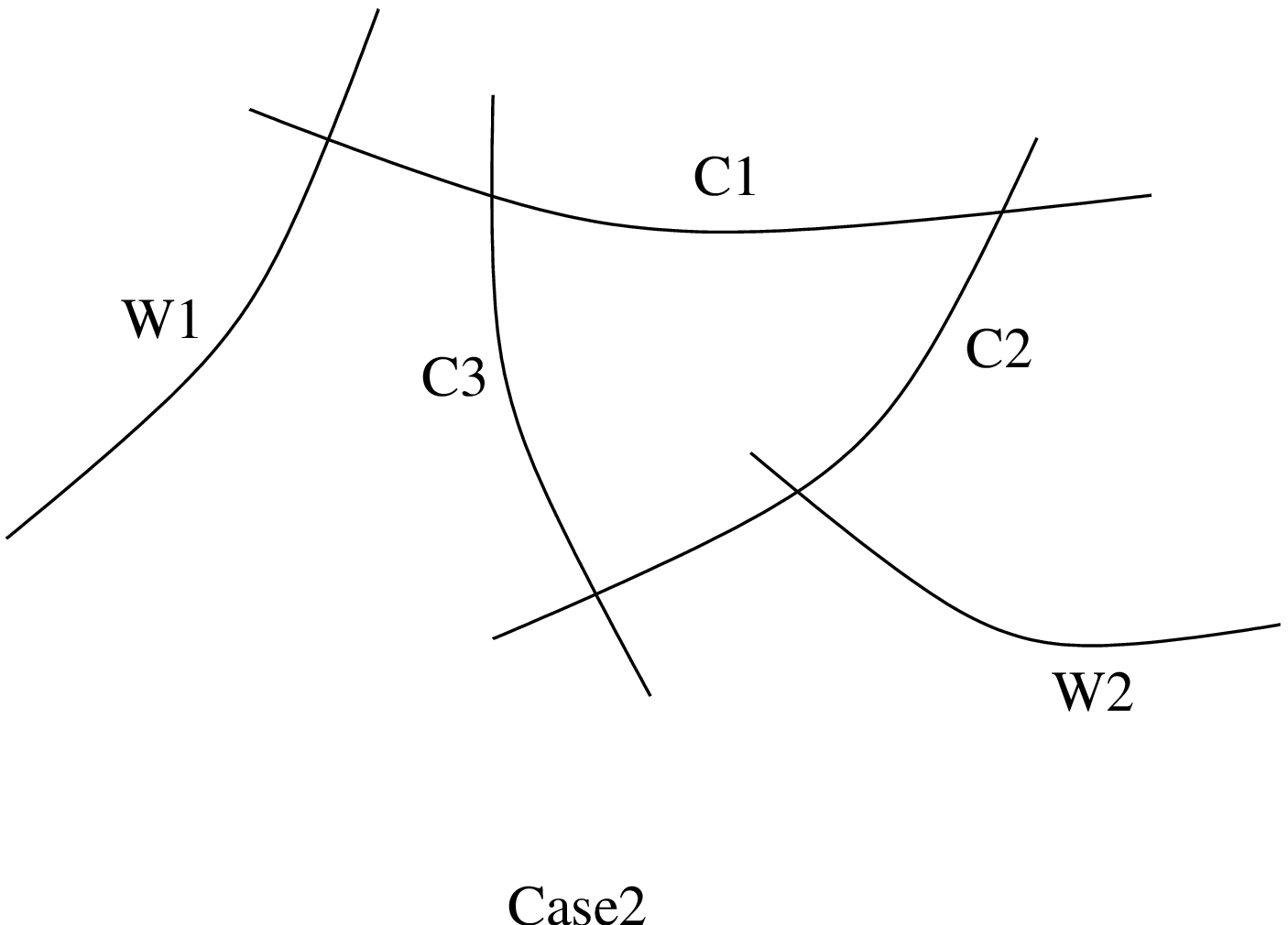}}

\end{center}
\end{enumerate}
We shall see that the numerics of $S_X$ is the same in all cases.

In case (a) there are two even sets of nodes: $\emptyset$ and $\{N\}$. If $\ev =\emptyset$ then
$b_1(\ev)=0$ and we find $2^{2g-2}$ components of $S_X$, all having multiplicity $2$,
supported on the quasistable curve obtained by ``blowing up " all nodes of $X$.

If $\ev =\{N\}$ then $b_1(\ev)=1$
and we find $2^{2g-1}$ reduced components of $S_X$
supported on the quasistable curve obtained by ``blowing up " all nodes of $X$ but $N$.

In case (b) there are again two even sets of nodes: $\emptyset$ and $D_{sing}$ (the $h$ nodes
of $D$). If $\ev =\emptyset$ then
$b_1(\ev)=0$ and, just as before, there are  $2^{2g-2}$ components of $S_X$, all having multiplicity $2$.

If $\ev =D_{sing}$ then $b_1(\ev)=1$
and  there are  $2^{2g-1}$ reduced components of $S_X$,
supported on the quasistable curve obtained by ``blowing up " all nodes of $X$, with the exception
 of  the nodes of $D$.
\end{example}
The previous example is special in the sense that the
numerical data of $S_X$ are independent of $X$ (so long as $b=1$). We
leave it to the reader  to check 
on other examples that
for $b\geq 2$ this is no longer true.
The picture will be made clear in the sequel (see Example \ref{es1} in the next section).

\begin{example}[Split curves]
\label{split}
A {\it split curve} $X$ is  defined to be a stable curve 
made of two irreducible components, $X=C_1\cup C_2$, with $C_i\cong {\mathbb {P}} ^1$; 
thus $C_1$ and $C_2$ meet in $g+1$ distinct points (the nodes of $X$),
and $b=g$.

Clearly, a set $\ev$  of nodes of $X$ is even if and only if its cardinality is even.
Denote by $d = \#\ev$, so that $b_1(\ev )=d-1$, unless $d=0$ in which case
$b_1(\emptyset )=0$.

Since $X$ has   $\binom{g+1}{d}$ even subsets of cardinality $d$, 
we see that, if $d$ is even and $0<d\leq g+1$, then  $S_X$ 
possesses $\binom{g+1}{d}$  irreducible components of multiplicity $2^{g-d+1}$;
for $d=0$, $S_X$ has a unique component of multiplicity $2^g$.
\end{example}
\begin{defi}
\label{mult} 
Let $S$ be a scheme of pure dimension zero. Denote by $L(S)$ the {\it multiplicity set} of $S$, that is,
set of
integers occurring as multiplicities of  components of  
$S$.
\end{defi}
For example, $S$ is reduced if and only if $L(S)=\{1\}$. 

We have, for
a split curve $X$ of genus $g$, 
$$
L(S_X) =\begin{cases}\{1,2^2,\dotsc,2^{g-3},2^{g-1}, 2^g\}\quad \text{ if $g$  is odd, }
\cr \{2,2^3,\dotsc ,2^{g-3},2^{g-1},2^g\}
\quad \text{ if   $g$  is even. }\cr    \end{cases}
$$
We have seen in Example \ref{ct} that curves of compact type are characterized by the fact that
the multiplicity set of their scheme of spin curves is equal to $\{1\}$.

What is remarkable is that the 
analogue holds for split curves (Corollary \ref{duec}).
This is a  consequence of a stronger result about the combinatorics of stable curves,
to which the next section will be devoted.

\section{Graph theoretic results and applications}
 \subsection{Preliminaries}
\label{cycles} 
Let $\g$ be a graph. We allow $\g$ to have loops and multiple edges,
namely: an edge can join a vertex to itself, and more than one edge can
join two vertices (these kinds of graphs are sometimes called 
multigraphs in the literature).
We also assume that $\g$ has no isolated vertices.

We recall some basic definitions and properties,
details can be found in  \cite{harary,diestel}, for example.

We  denote by $E(\g)$ and $V(\g)$
respectively the sets of edges and vertices of $\g$. Let $\delta_{\g}=\#
E(\g)$, $\nu_{\g}=\#V(\g)$ and let $c(\g)$ be the number of connected
components of $\g$. The first Betti number of $\g$ is
$b_1(\g)=\delta_{\g}-\nu_{\g}+c(\g)$. 
We shall say that a vertex (or an edge) $a$ is a {\it separating vertex} (a
{\it separating edge}) if $c(\g\smallsetminus\{a\})>c(\g)$.

A subgraph $\g'\subseteq\g$ is a graph $\g'$ such that $V(\g')\subseteq V(\g)$
and $E(\g')\subseteq E(\g)$. Given two graphs $\g_1,\g_2$ their union
is the graph $\g_1\cup\g_2$ such that
$V(\g_1\cup\g_2)=V(\g_1)\cup V(\g_2)$ and $E(\g_1\cup\g_2)=E(\g_1)\cup
E(\g_2)$. 

The valency of a vertex is the number of edges ending in
that vertex (a loop counting 2 in the valency). 
Clearly the sum of all valencies equals $2\delta_{\g}$. 

 A \emph{circuit} of $\g$ is a subgraph of $\g$ that has all valencies
equal  to 2.

We consider a vector space associated to $\g$, its \emph{cycle space}
$\mathcal{C}_{\g}$. This is a vector space over
$\mathbb{F}_2=\Z/2\Z$, of dimension $b_1(\g)$. A 1-chain (respectively,
0-chain) of $\g$ is a
formal linear combination of edges (respectively, vertices) of $\g$ with
coefficients in $\mathbb{F}_2$. 

Any 1-chain $\ev=N_1+\cdots+N_m$ can be viewed 
as a subset
$\{ N_1,\dotsc,N_m\}$ of $E(\g)$, or as a subgraph of $\g$, namely the 
smallest subgraph of $\g$ containing
$N_1,\dotsc,N_m$. 
To avoid a heavy notation, we will
not make any distinction among these different interpretations, 
 denoting by the same letter, $\ev$, the 1-chain, the set of
edges and the induced graph.

The boundary operator $\partial$ sends
1-chains to 0-chains in the usual way: $\partial$ is linear and for an
edge $N$, $\partial(N)=a+b$, where $a$ and $b$ are the vertices of
$N$. The cycle space $\mathcal{C}_{\g}$ is the kernel of the boundary
operator (using the standard notation, 
$\mathcal{C}_{\g}:=Z_1(\g, \mathbb{F}_2)$).
  An element $\ev\in\mathcal{C}_{\g}$ is called a  \emph{cyclic
  set}  or a  \emph{cyclic
  subgraph} of $\g$. 
\begin{claim}
A subgraph $\ev\subseteq\g$ is  cyclic  if and only if $\ev$
has all even valencies.
\end{claim}
\begin{proof}
Consider 1-chains and 0-chains with coefficients in $\Z$ instead
of $\mathbb{F}_2$, and consider the 0-chain $\partial(\ev)$.
The coefficient of any vertex in
$\partial(\ev)$ is exactly its valency in $\ev$. Hence the
statement follows.
\end{proof}
In particular, a circuit is a cyclic set, and every cyclic
set can be decomposed as an edge-disjoint union of circuits.

A graph $\g$ having all even valencies, or equivalently such that
$E(\g)\in\mathcal{C}_{\g}$, is called an \emph{eulerian} graph. Such a
graph is characterized by the existence of a closed walk passing
exactly once through every edge. 

Let now $X$ be a  stable curve and $\g = \g_X$ be its dual graph.
Let $\ev\subset X_{sing}$ 
be a subset of nodes and (with the usual abuse of notation)
$\ev \subset \g$ the corresponding subgraph.
\begin{remark}
To say that $\ev$ is even in the sense of Definition~\ref{even} is exactly the same as saying that $\ev \in
\mathcal{C}_{\g_X}$, or that $\ev$ is cyclic, 
in the graph theoretic language. In particular, we see that on $X$ there are exactly $2^{b_1(\g_X)}$ even subsets of
nodes.
\end{remark}

 \subsection{Relating to spin curves}
As we saw in the previous section, the geometry of $S_X$ is ruled by the even sets of nodes of the stable curve $X$, that is, by the
cyclic subgraphs of
$\g _X$. 
\begin{defi}
\label{betti}
The \emph {set  of cyclic Betti numbers} of a graph $\g$ 
is  $B_{\g}:=\{b_1(\ev)\,|\,\ev\in\mathcal{C}_{\g}\}$.
\end{defi}
The relevance of such a definition comes from Proposition
\ref{num}. In fact,  let  
$b=b_1(\Gamma _X)$,
we have (see Definition \ref{mult})
$$
L(S_X)=\{2^{b-n}\,|\, n \in B_{\g _X}\}.
$$
For example, the analysis of Example \ref{ct} shows that
$X$ is of compact type if and only $B_{\g _X}=\{ 0\}$. 

\paragraph{ Elementary properties of  $B_{\g}$.}
\begin{enumerate}[(P1)]
\item \label{a}
\emph{For all $m\in B_{\g}$, we have $m\leq b_1(\g)$}.

In fact, for every subgraph $\ev\subset\g$, clearly $b_1(\ev)\leq
b_1(\g)$.
\item \label{b}
$0\in B_{\g}$. 

The vector space $\mathcal{C}_{\g}$ 
contains the zero element.
\item \label{c}
\emph{$B_{\g}=\{0\}$ if and only if $\ \g$ is a tree, if and only if
  $\,1\not\in B_{\g}$}.

 A circuit $\gamma$ always has $b_1(\gamma)=1$, and 
$b_1(\g)\geq 1$ if and only if $\g$ contains a circuit, if and only if 
$1\in B_{\g}$. 
\item \label{d}
\emph{If $\g'\subset\g$, then $\mathcal{C}_{\g'}$ is a subspace of 
$\mathcal{C}_{\g}$ and $B_{\g'}\subseteq B_{\g}$}.

 The vector space of 1-chains of $\g'$ is a subspace
of the 1-chains of $\g$.  
\item \label{e}
\emph{If $\g=\g_1\cup\g_2$, $E(\g_1)\cap E(\g_2)=\emptyset$ and  
$\#V(\g_1)\cap V(\g_2)\leq 1$, then $\mathcal{C}_{\g}=
\mathcal{C}_{\g_1}\oplus\mathcal{C}_{\g_2}$ and $B_{\g}=B_{\g_1}+B_{\g_2}=
\{n_1+n_2\,|\,n_1\in B_{\g_1},n_2\in B_{\g_2}\}$.}

Since $\g=\g_1\cup\g_2$ and $E(\g_1)\cap E(\g_2)=\emptyset$, we
have $\mathcal{C}_{\g}= \mathcal{C}_{\g_1}\oplus\mathcal{C}_{\g_2}$.
Since $\#V(\g_1)\cap V(\g_2)\leq 1$,
given two cyclic sets $\ev_1\in\mathcal{C}_{\g_1}$ and
$\ev_2\in\mathcal{C}_{\g_2}$,  either
they are disjoint or they have one common vertex. In both cases we
have $b_1(\ev_1\cup\ev_2)=b_1(\ev_1)+b_1(\ev_2)$. 
\item \label{f}
\emph{ A cyclic set does not contain any separating edge.}
 
An edge $N$ is contained in some cyclic set if and only if
it is not a separating edge: in fact, if $a$ and $b$ are the vertices
of $N$, $N$ is not a separating edge if and only if there is a path in
$\g\smallsetminus\{N\}$ joining $a$ and $b$. This path together with
$N$ gives a cyclic set in $\g$. 
\item \label{g}
\emph{$b_1(\g)\in B_{\g}$ if and only if the set
  $\,E(\g)\smallsetminus\{\,\text{separating edges}\,\}$ is a cyclic set.}

Set $r=\#\{\,\text{separating edges}\,\}$ and
$\g'=\g\smallsetminus\{\,\text{separating edges}\,\}$. Then we have 
$\delta_{\g'}=\delta_{\g}-r$, 
$\nu_{\g'}=\nu_{\g}$ and $c(\g')=c(\g)+r$, so $b_1(\g')=b_1(\g)$.
Hence, if $\g'$ is a cyclic set, then $b_1(\g)\in B_{\g}$.
 
Suppose now that $\ev$ is a cyclic set strictly smaller than 
$\,E(\g)\smallsetminus\{\,\text{separating edges}\,\}$: it suffices to show that
$b_1(\ev)<b_1(\g)$. In fact,
there exists a non-separating edge $N$ such that 
 $\ev$ is contained in the
subgraph $\g'=\g\smallsetminus\{N\}$; we have $\delta_{\g'}=\delta_{\g}-1$,
$\nu_{\g'}=\nu_{\g}$ and $c(\g')=c(\g)$, so $b_1(\ev)\leq
b_1(\g')=b_1(\g)-1$. 
\item \label{h}
\emph{A graph $\g$ is eulerian if and only if it has no
separating edges and $b_1(\g)\in B_{\g}$.}

 This follows from properties (P\ref{f}) and (P\ref{g}).
\end{enumerate}

\paragraph{Operations on $\g$ that fix $\mathcal{C}_{\g}$ and $B_{\g}$:}
\begin{enumerate}[1.]
\item elimination of a vertex of valency 1: we contract the unique edge
containing the vertex;
\item
elimination of a vertex of valency 2, not allowed on the vertex of a
loop: the two edges ending in the vertex are merged
in a unique edge;
\item
elimination of a separating edge: the edge is contracted and gives
rise to a separating vertex.
\end{enumerate}
\begin{lemma}
\label{operazioni}
If $\g'$ is obtained from $\g$ by any sequence of operations of type 1,2
and 3, we have $b_1(\g')=b_1(\g)$,
$\mathcal{C}_{\g'}\simeq\mathcal{C}_{\g}$ 
and $B_{\g'}=B_{\g}$. 
\end{lemma}
\begin{proof} Elementary and easy.
\end{proof}

\begin{defi}
A graph $\g$ is \emph{superstable} if  all valencies of $\g$ are at least 3, except possibly for
the vertex of a loop.
\end{defi}
For any graph $\g$, there exists a unique superstable graph $\g^s$ obtained
from $\g$ with a sequence of operations of type 1 and 2, so
that $\mathcal{C}_{\g^s}\simeq\mathcal{C}_{\g}$ and 
$B_{\g^s}=B_{\g}$. 

To explain the choice of the name ``superstable", 
pick a nodal 
connected curve
$Z$, all of whose components have geometric genus zero and such that $\g$ is the dual graph of $Z$.
The fact that such a $Z$ may not be uniquely determined by $\g$ is irrelevant in the following discussion.
The point is that  $\g$ is  superstable if and only if the curve $Z$ is stable (in the sense of Deligne and Mumford).
Notice moreover that the three above operations on $\g$  
correspond to operations on $Z$: the first contracts a smooth tail, the second contracts an exceptional 
 component (in the sense of Definition \ref{stable}), the third corresponds to smoothing a separating node. 
Hence operations  1 and 2 correspond to ``stabilizing" operations on the curve;  in other words,
the   graph $\g^s$ defined above is the dual graph of the stable model of $Z$ (in the sense of Definition
\ref{stable}).  

\begin{example}[superstable graphs with $b_1(\g)=1$ or $2$]
\label{es1}
Let us 
only consider  graphs free from separating edges, leaving it to the reader to list the
remaining ones. It is immediate to see that if $b_1(\g)=1$, 
$\g$ must be a loop, and $B_{\g}=\{0,1\}$ (this  clarifies Example \ref{b1} in the previous
section).

If $b_1(\g)=2$, either $\g$ has two loops as connected components, or
it is connected and has all valencies at least 3. 
In this last case, it is easy to see
that the only possibilities for $(\delta_{\g},\nu_{\g})$ are $(2,1)$
and $(3,2)$, and that there is only one possible graph for each pair. We conclude that 
the superstable graphs with $b_1(\g)=2$ and no separating edges are:
\begin{center}  

 \scalebox{0.4}{\includegraphics{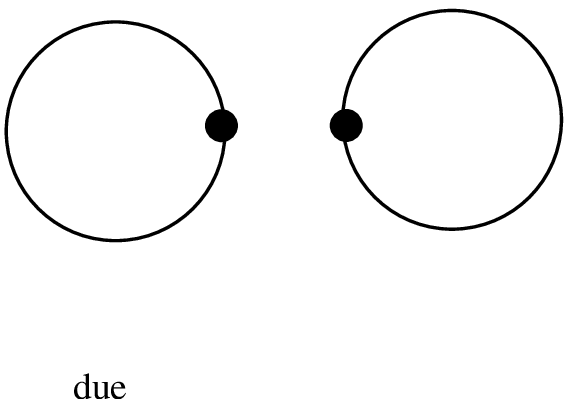}}
 \hspace{80pt}\scalebox{0.4}{\includegraphics{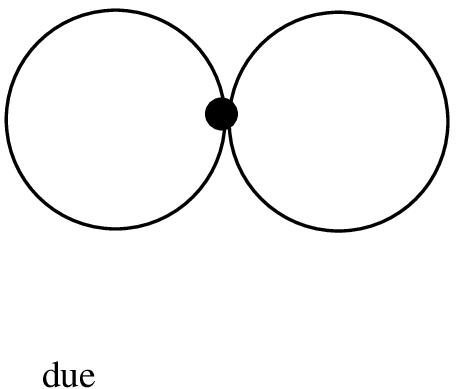}}
 \hspace{80pt}\scalebox{0.4}{\includegraphics{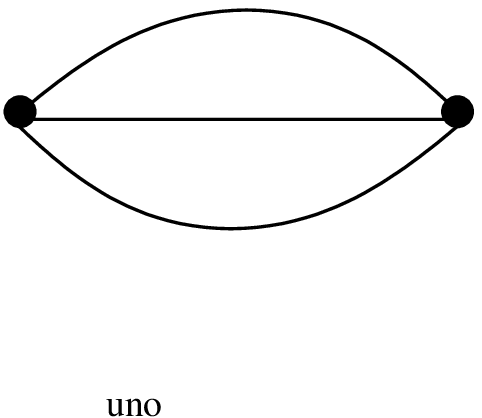}}

\end{center}
\end{example}
\noindent Now we introduce split graphs (compare with Example \ref{split}):
\begin{defi}
A graph is \emph{split} if it is connected, has two vertices and
no loops.
\end{defi}
If $\g$ is split and $b_1(\g)\geq 2$, then $\g$ is superstable.

Let $\g$ be split and consider a subgraph $\ev\subset\g$ 
containing $r$ edges. Then $b_1(\ev)=r-1$, and
$\ev$ is a cyclic set if and only if $r$ is even. 
Hence \[B_{\g}=\{0,k\,|\,k\text{ is odd and }k\leq
b_1(\g)\}. \]

\subsection{Combinatorial results} 

We now prove two structural characterizations of graphs, using their cyclic Betti numbers
(Theorem~\ref{due} and Theorem~\ref{tre}). 

Parenthetically, we mention that  Theorem~\ref{due} and Corollary~\ref{duec} 
  generalize (and clarify) a crucial  step in  \cite{cs}
(Theorem 3.4.1); the main result of that paper is the
fact that
 a general
 canonical curve is uniquely determined by the configuration of
hyperplanes cutting theta characteristics on it. 
Applying stable reduction, such a classical, concrete problem is solved using the moduli theory
of stable curves.
The combinatorial analysis is  used as a bridge between the projective and the
abstract set-up.

\begin{thm}
\label{due}
Let $\g$ be a superstable graph such that $2\not\in B_{\g}$. Then
either $\g$ is split, or $b_1(\g)=1$ and $\g$ is a loop, or
$b_1(\g)=3$ and $\g$ is the tetrahedron:

 \begin{center}
  
 \scalebox{0.4}{\includegraphics{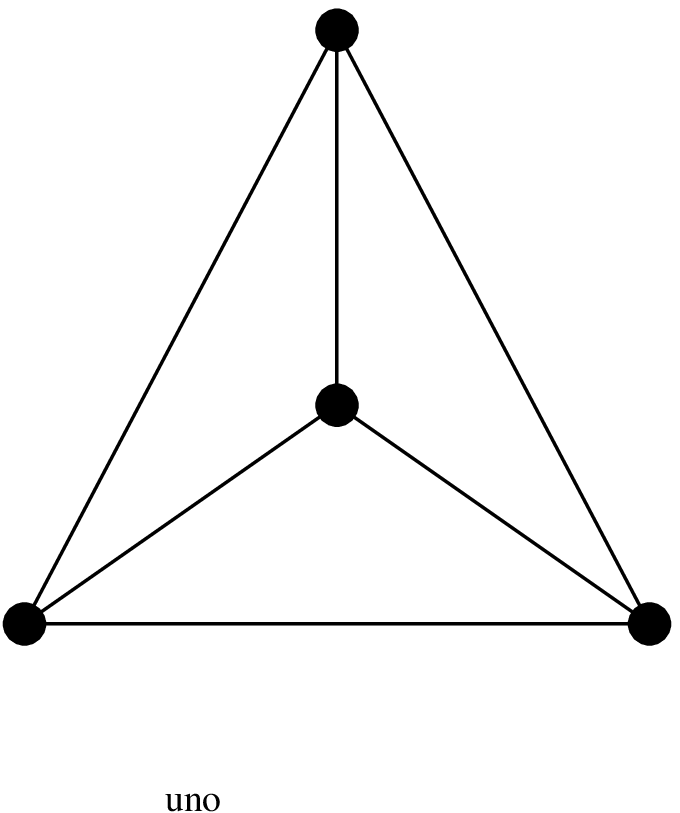}}

 \end{center}

\end{thm}
\begin{proof}
We remark first of all that
$\g$ is connected and does not have separating vertices. Otherwise,
there would be two subgraphs $\g_1$, $\g_2$ such that
$\g=\g_1\cup\g_2$, $E(\g_1)\cap E(\g_2)=\emptyset$ and  
$\#V(\g_1)\cap V(\g_2)\leq 1$. Since $\g$ is superstable, $\g _i$ cannot be a tree,
for
$i=1,2$. By property (P\ref{c}),
$b_1(\g_i)\geq 1$ and 
$B_{\g_1}$ and $B_{\g_2}$ contain $1$. 
Property (P\ref{e}) implies 
$2\in B_{\g}=B_{\g_1}+B_{\g_2}$, a contradiction.

From the fact that $\g$ has no separating vertices, we deduce that
$\g$ has no separating edges. If it did, let $\g '$ be the superstable graph 
obtained by contracting all separating edges; since every such  contraction
generates a separating vertex,
$\g '$ would possess separating vertices.
On the other hand, by Lemma \ref{operazioni}, $B_{\g}=B_{\g'}$, therefore $\g'$ satisfies the Theorem's
assumptions,
hence, by what we proved above, $\g'$ has no separating vertices; a contradiction.

We proceed by induction on $b_1(\g)$. If $b_1(\g)=1$ or $2$, the 
statement follows from Example~\ref{es1}.

Let's suppose  $b_1(\g)\geq 3$. We choose an edge $N\in E(\g)$ and
consider the subgraph
$\widetilde{\g}=\g\smallsetminus\{N\}$. $\widetilde{\g}$ could have one or
two vertices of valency 2; we eliminate those vertices (operation 2) and
obtain a connected, superstable graph $\g'$ such that
$b_1(\g')=b_1(\widetilde{\g})=b_1(\g)-1$. 
By Lemma~\ref{operazioni} and property (P\ref{d}) we have
$B_{\g'}=B_{\widetilde{\g}}\subseteq B_{\g}$, so $2\not\in
B_{\g'}$: the induction hypothesis implies that  either $\g'$ is split, or
$b_1(\g)=4$ and $\g'$ is the tetrahedron.

Suppose first that $\g'$ is split. If $\g$ is not split, there are
three possibilities for $\g$:
\begin{center}
  
\scalebox{0.25}{\includegraphics{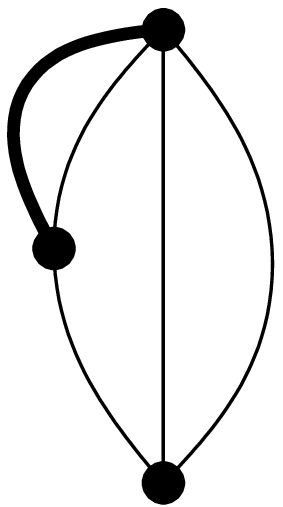}}\hspace{90pt}\scalebox{0.25}{\includegraphics{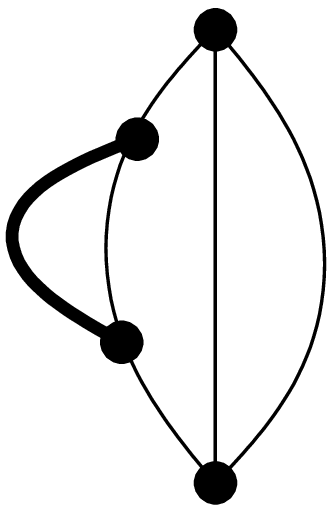}}
\hspace{90pt}\scalebox{0.25}{\includegraphics{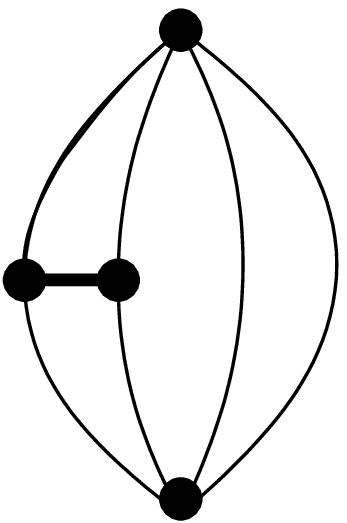}}

\end{center}
The thickened edge in the picture is $N$.
Recall that $b_1(\g')\geq 2$. 

In the first two cases we obtain a  contradiction
because there clearly is a
cyclic set $\ev\subset\g$
such that $b_1(\ev)=2$:
\begin{center}
  
\scalebox{0.25}{\includegraphics{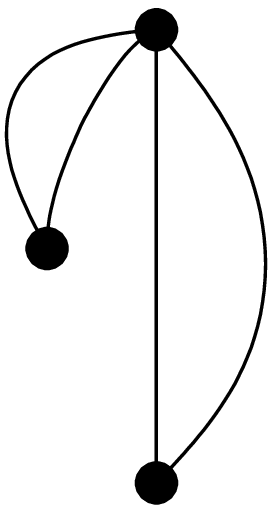}}\hspace{90pt}\scalebox{0.25}{\includegraphics{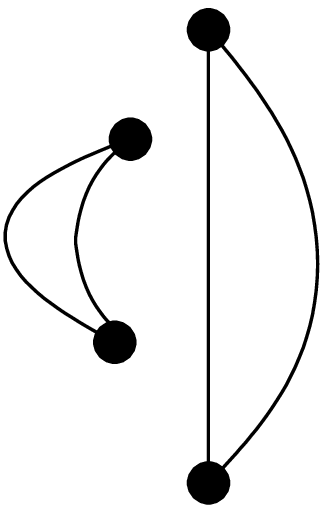}} 

\end{center}
In the third case,   if $b_1(\g')\geq 3$, we can find again a cyclic
set $\ev\subset\g$ with $b_1(\ev)=2$
(picture below, on the left), while if $b_1(\g')= 2$ we get the
tetrahedron (picture below, on the right):
\begin{center}
  
\scalebox{0.25}{\includegraphics{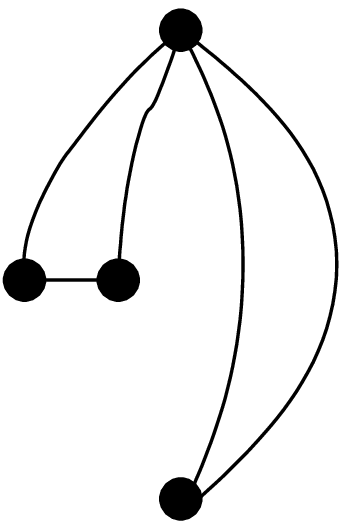}}\hspace{120pt}\scalebox{0.25}{\includegraphics{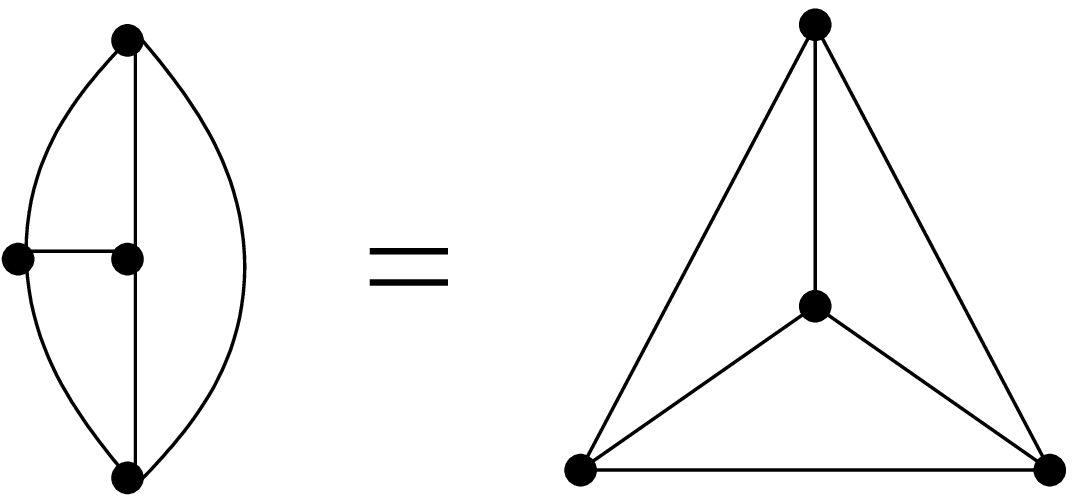}}

\end{center}

Suppose now that $\g'$ is the tetrahedron: then in any case we get a
cyclic set $\ev\subset\g$ such that $b_1(\ev)=2$:
 \begin{center}
  
$\g:$\hspace{30pt}
\scalebox{0.2}{\includegraphics{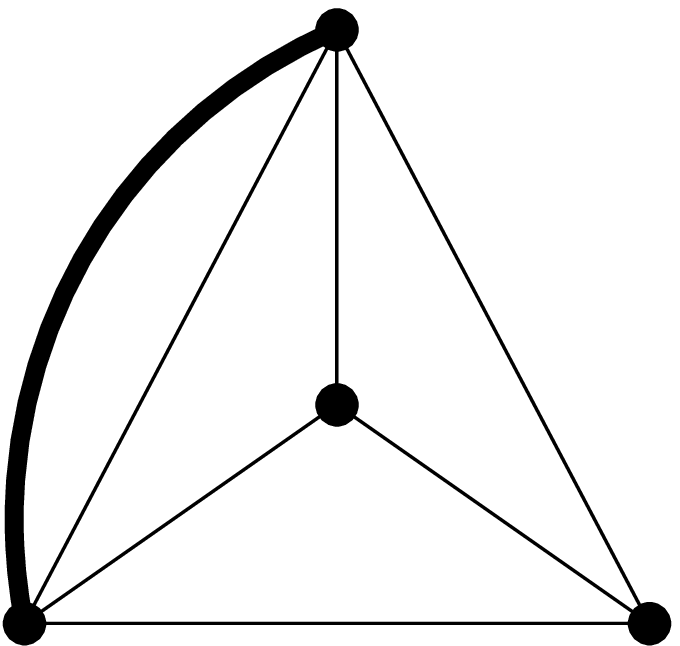}}\hspace{70pt}\scalebox{0.2}{\includegraphics{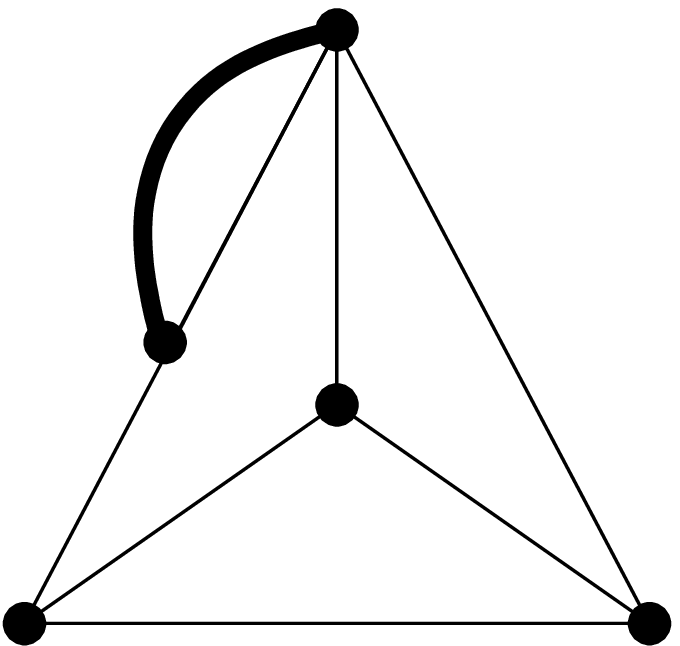}}
\hspace{70pt}\scalebox{0.2}{\includegraphics{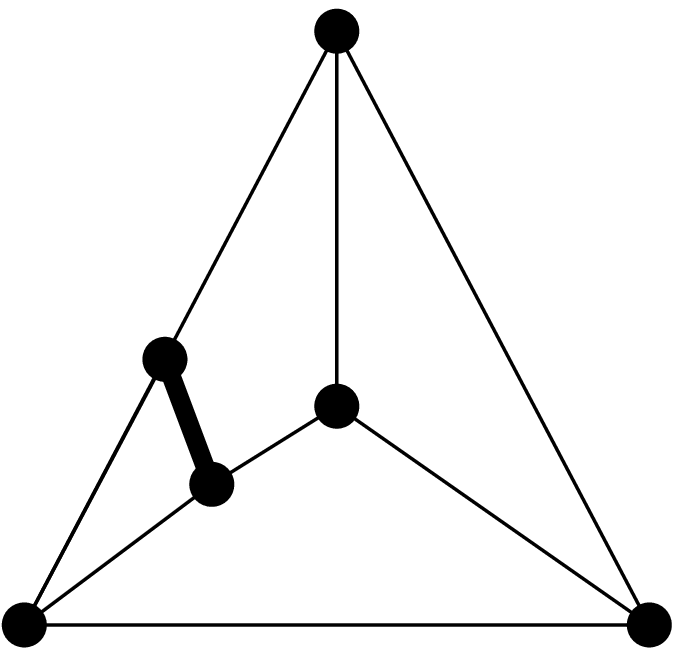}}

\vspace{20pt}

$\Delta:$\hspace{35pt}\scalebox{0.2}{\includegraphics{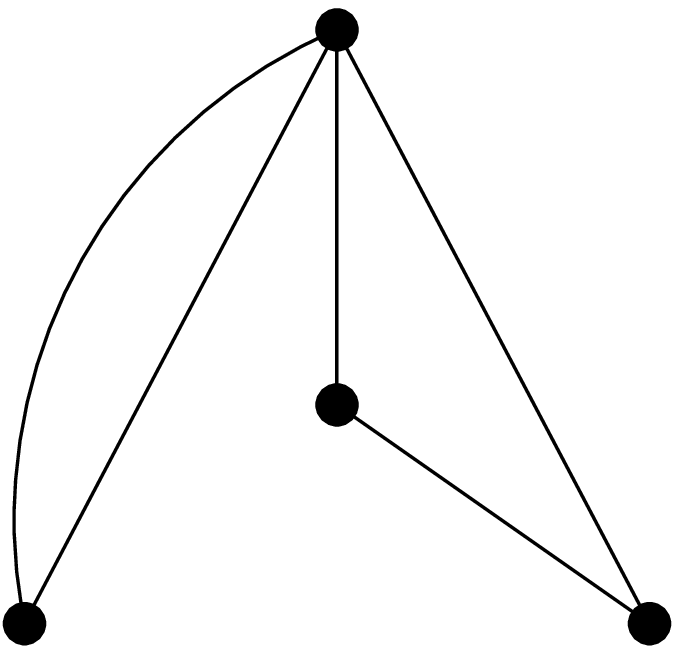}}\hspace{75pt}
\scalebox{0.2}{\includegraphics{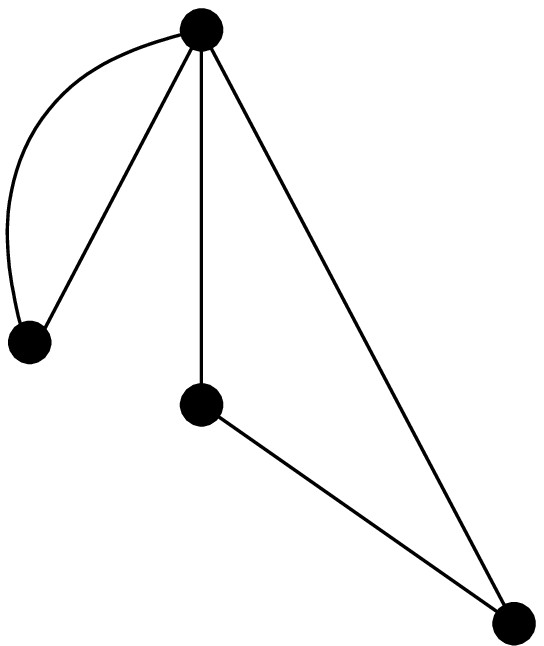}}\hspace{70pt}\scalebox{0.2}{\includegraphics{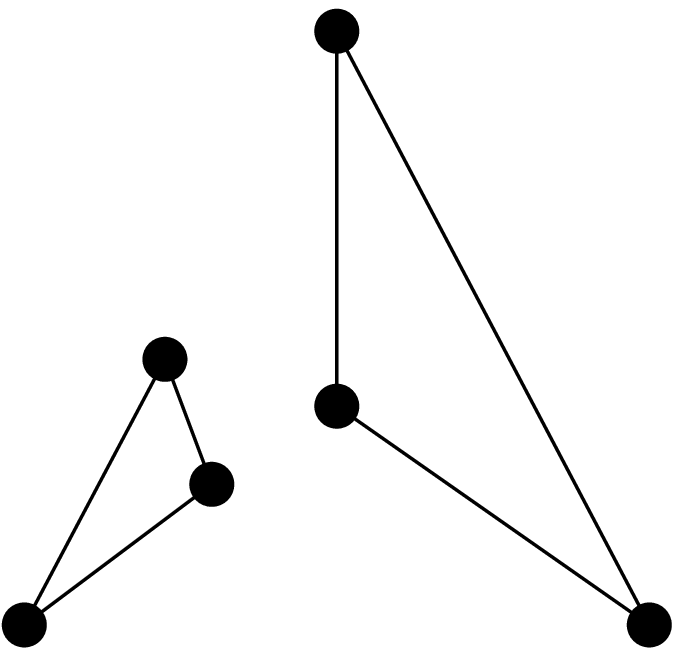}}

 \end{center}
\end{proof}
Applying this to  stable  curves we obtain:
\begin{cor}
\label{duec}
Let $X$ be a stable curve of genus $g$. Suppose that $S_X$ has a  component of multiplicity $2^g$ and no component of
multiplicity $2^{g-2}$. Then either $X$ is a split curve or $X$ has genus $3$ and it is the polygonal curve.
\end{cor}
\noindent
(Recall that the polygonal curve of genus $3$ is the one whose dual graph is a tetrahedron.)
\begin{proof}
The fact that $2^g\in S_X$ implies, by Proposition \ref{num}, that $b_1(\g _X)=g$ and hence every irreducible
component of
$X$ has geometric genus $0$. Since $X$ is stable, $\g_X$ must have all valencies at least equal to $3$.
To say that $2^{g-2}\not\in S_X$ is the same as saying that
$2\not\in B_{\g_X}$. Therefore all
the assumptions of Theorem \ref{due} are satisfied and $X$ is a split curve or its dual graph is the tetrahedron.
\end{proof}

\begin{thm}
\label{tre}
Let $\g$ be a superstable graph. Suppose  that $3\not\in B_{\g}$
and that there exists $m\in B_{\g}$ such that $m>3$. Then $b_1(\g)=4$ and
$\g$ is the ``fat-triangle'':
\begin{center}
  
  \scalebox{0.3}{\includegraphics{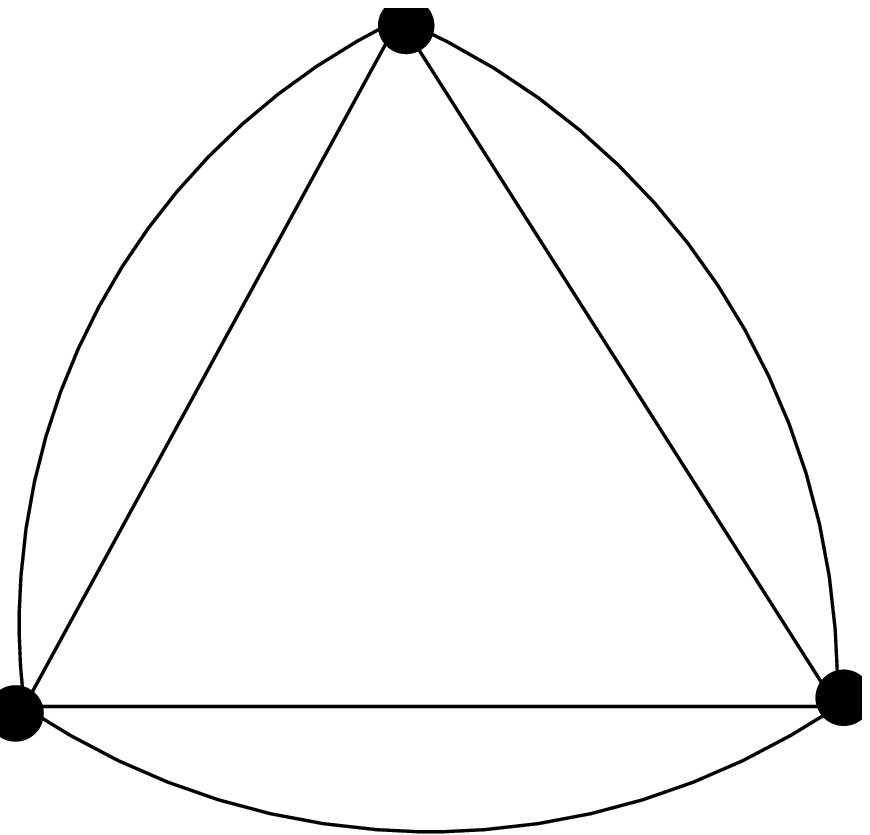}}

\end{center}
\end{thm}
\begin{proof}
Remark that $b_1(\g)\geq m\geq 4$, so by Theorem~\ref{due} we know
that $2\in B_{\g}$. Moreover,
$\g$ is connected and does not have separating vertices. 
Otherwise,
there would be two subgraphs $\g_1$, $\g_2$ such that
$\g=\g_1\cup\g_2$, $E(\g_1)\cap E(\g_2)=\emptyset$ and  
$\#V(\g_1)\cap V(\g_2)\leq 1$. Since $\g$ is superstable, we have
$b_1(\g_1)\geq 1$ and $b_1(\g_2)\geq 1$. 
Then  by properties (P\ref{b}) and (P\ref{c})
$B_{\g_1}$ and $B_{\g_2}$ contain 0 and 1, but since
$B_{\g}=B_{\g_1}+B_{\g_2}$ and $3\not\in B_{\g}$, 
we must have that $2,3\not\in B_{\g_i}$ for $i=1,2$. We deduce from
Theorem~\ref{due} that $B_{\g_i}=\{0,1\}$, so $B_{\g}=\{0,1,2\}$,
a contradiction.

Just as we did in the proof of Theorem \ref{due}, we deduce that $\g$ has no separating edges.

We continue by induction on $b_1(\g)$. When $b_1(\g)=4$, we have
$B_{\g}=\{0,1,2,4\}$. Since $4=b_1(\g)\in B_{\g}$, $\g$ is eulerian 
(property (P\ref{h}));  
in particular all valencies are at least 4.
Thus the possibilities for the number of edges 
and vertices $(\delta_{\g},\nu_{\g})$ are only $(4,1)$, $(5,2)$ and
$(6,3)$. Now, the only possible $\g$ having $b_1(\g)=4$, even valencies and
without separating vertices is the fat-triangle.

Suppose now $b_1(\g)>4$.
\begin{claim}
There exists an $m\in B_{\g}$
such that $3<m<b_1(\g)$. 
\end{claim}
Let's see how the claim implies the statement. Let $\ev\subset\g$ be
a cyclic set such that $b_1(\ev)=m$. Since $m<b_1(\g)$, there exists an
edge $N\in E(\g)\smallsetminus E(\ev)$. Consider the superstable 
graph $\g'$ obtained from $\g$ eliminating $N$, as in the proof of
Theorem~\ref{due}.  Then
$B_{\g'}\subseteq B_{\g}$, so
 $3\not\in B(\g')$, but clearly $\ev\subseteq\g'$, so $m\in B(\g')$. 
Since $b_1(\g')=b_1(\g)-1$, by the induction hypothesis we deduce that
$b_1(\g)=5$ and
$\g'$ is the fat-triangle. Then there are five possibilities for $\g$:

\begin{center}
  
\scalebox{0.15}{\includegraphics{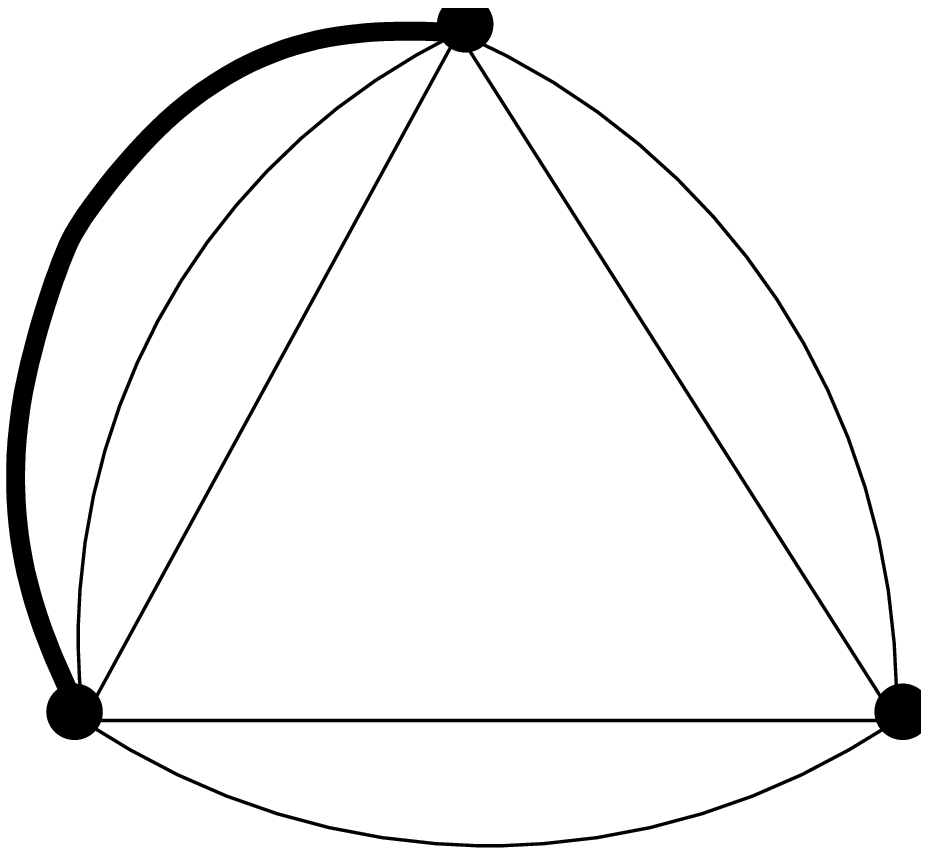}}
\hspace{35pt}\scalebox{0.15}{\includegraphics{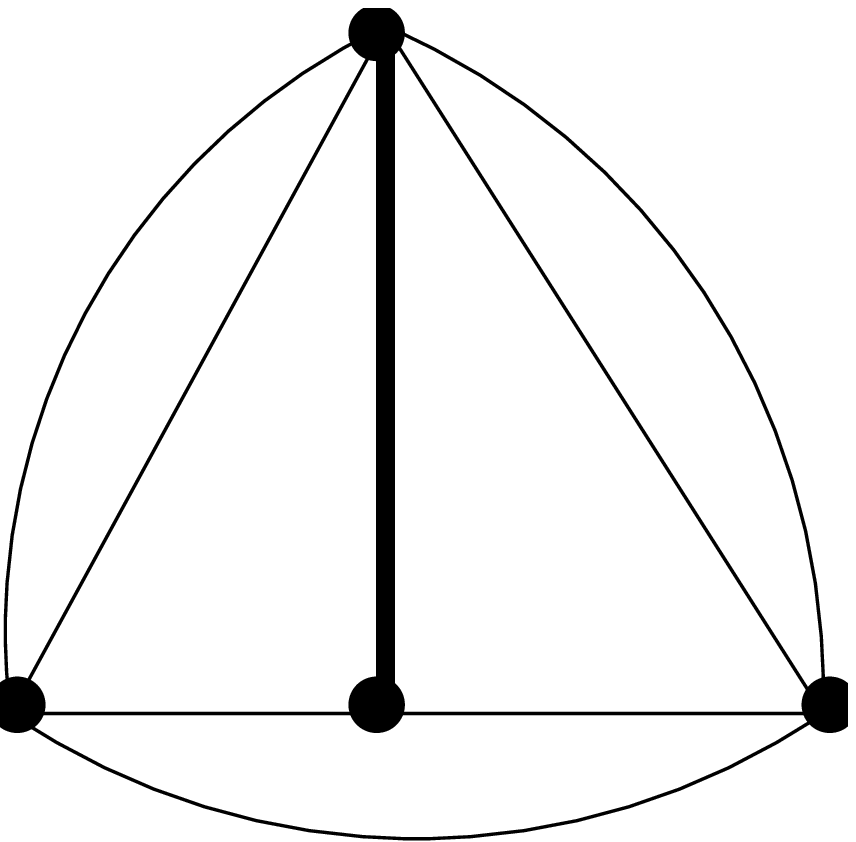}}\hspace{35pt}\scalebox{0.15}{\includegraphics{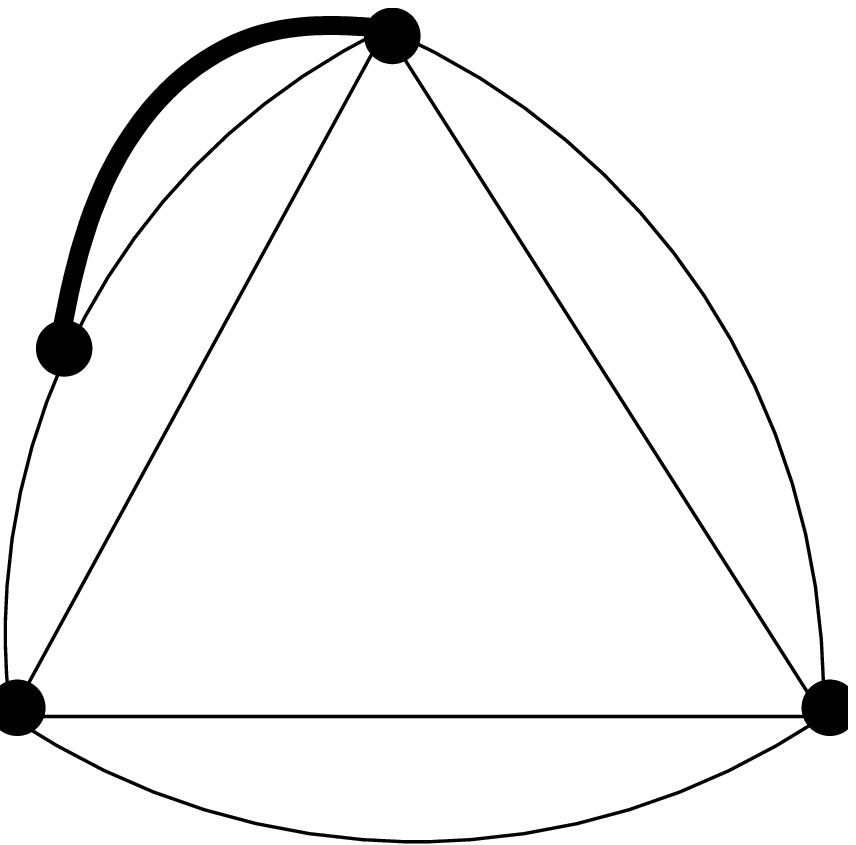}}
 \hspace{35pt}\scalebox{0.15}{\includegraphics{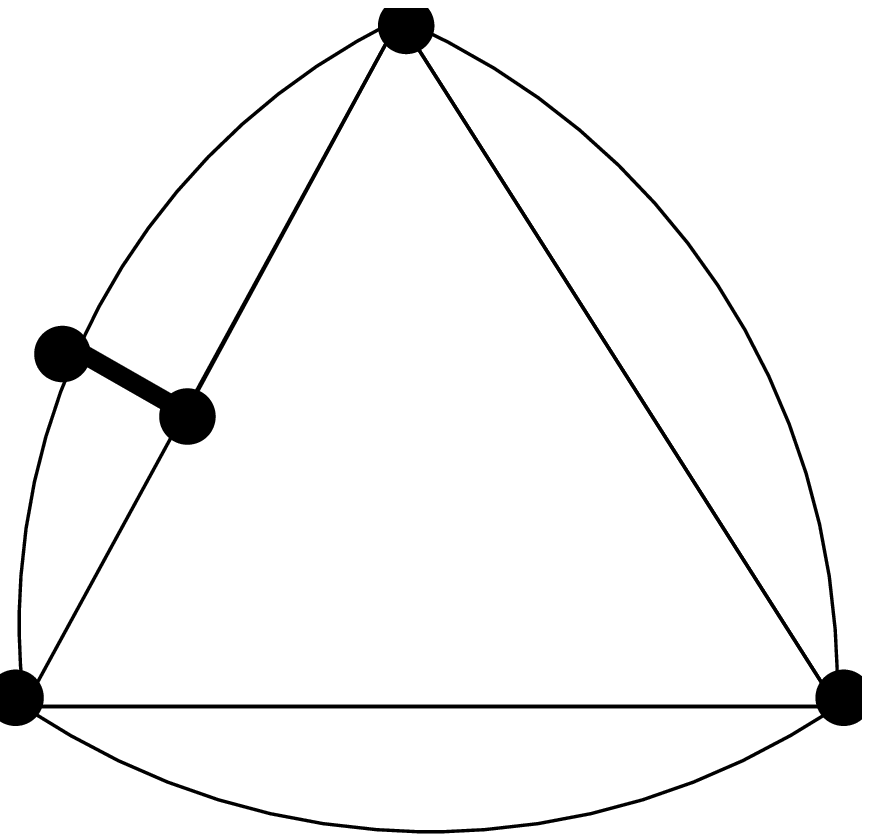}}\hspace{35pt}\scalebox{0.15}{\includegraphics{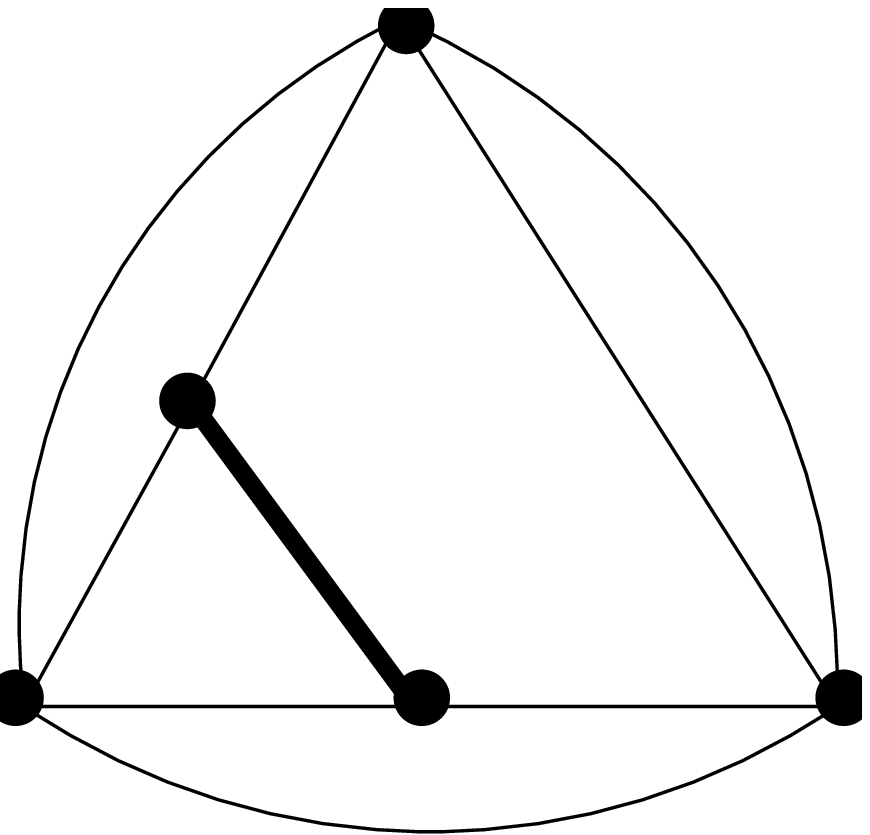}}

\end{center}

\noindent The thickened edge in the picture is $N$. In all these cases, there is a cyclic set $\ev$  with
$b_1(\ev)=3$, against the hypotheses:
\begin{center}
  
\scalebox{0.15}{\includegraphics{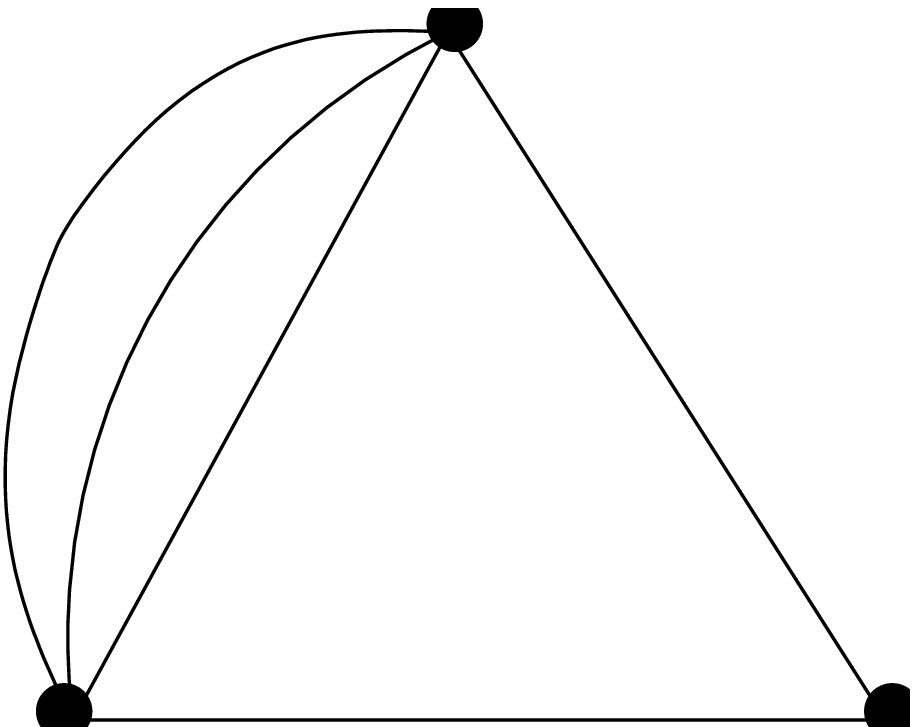}}
\hspace{35pt}\scalebox{0.15}{\includegraphics{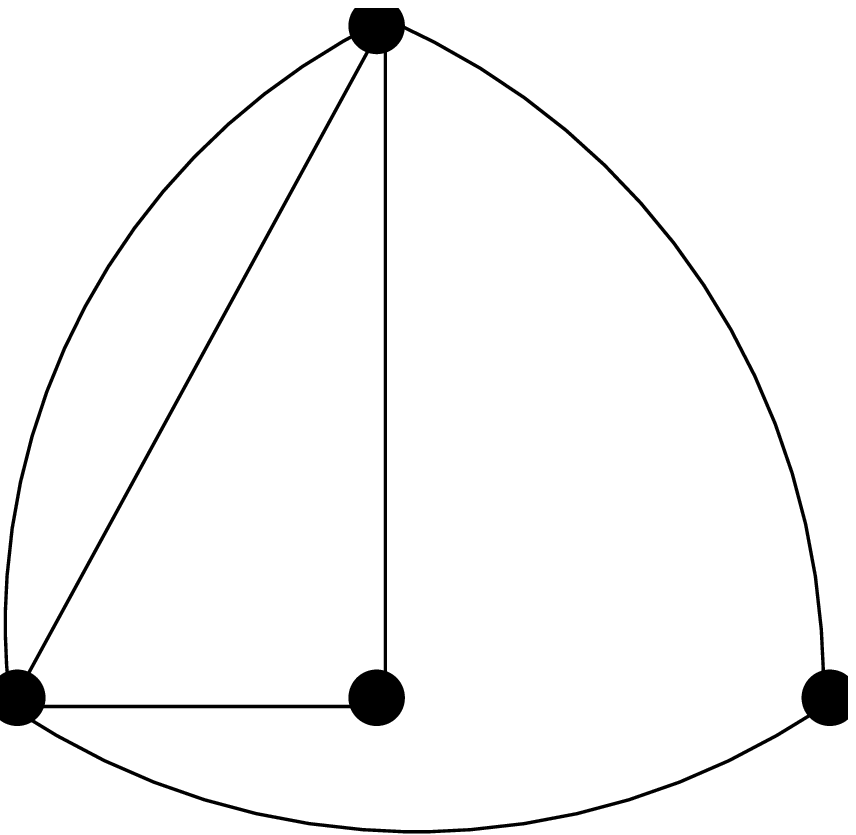}}\hspace{35pt}\scalebox{0.15}{\includegraphics{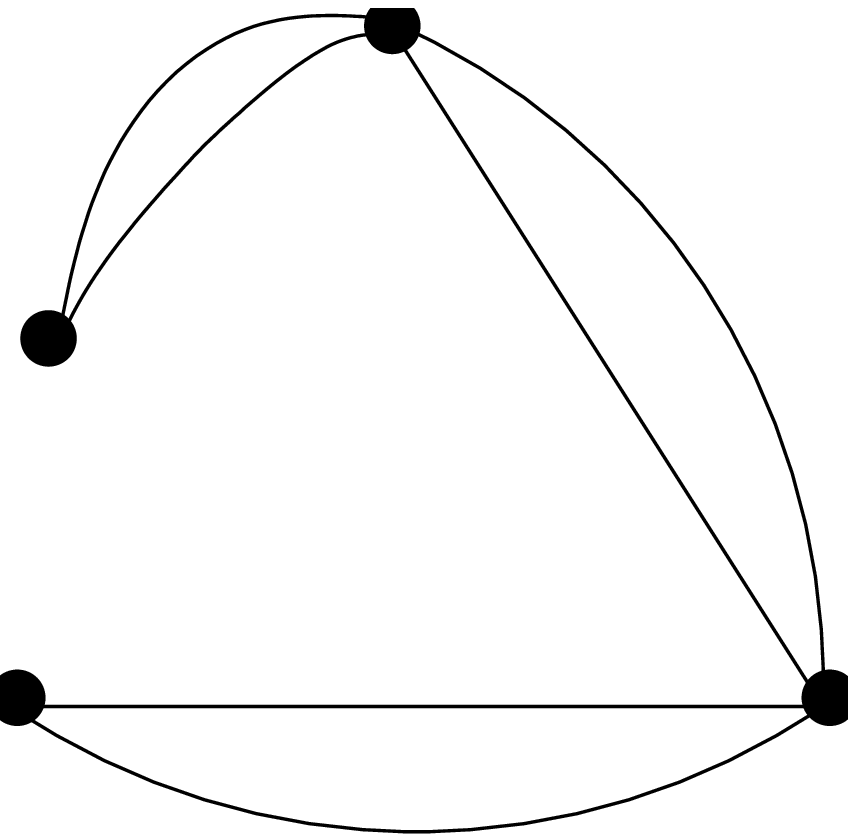}}
 \hspace{35pt}\scalebox{0.15}{\includegraphics{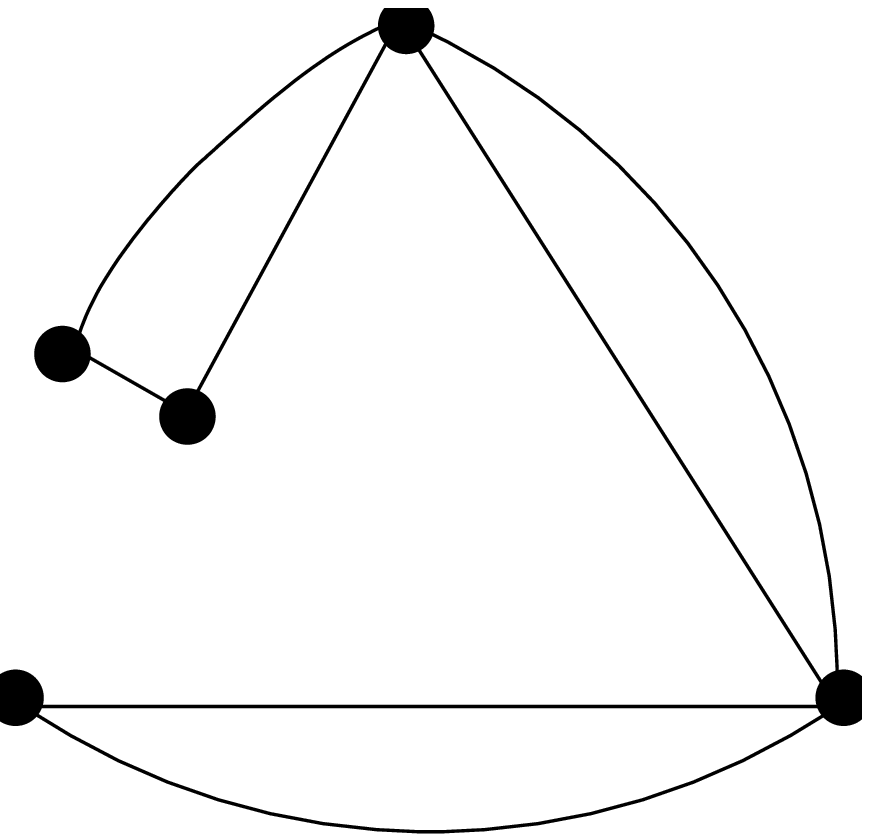}}\hspace{35pt}\scalebox{0.15}{\includegraphics{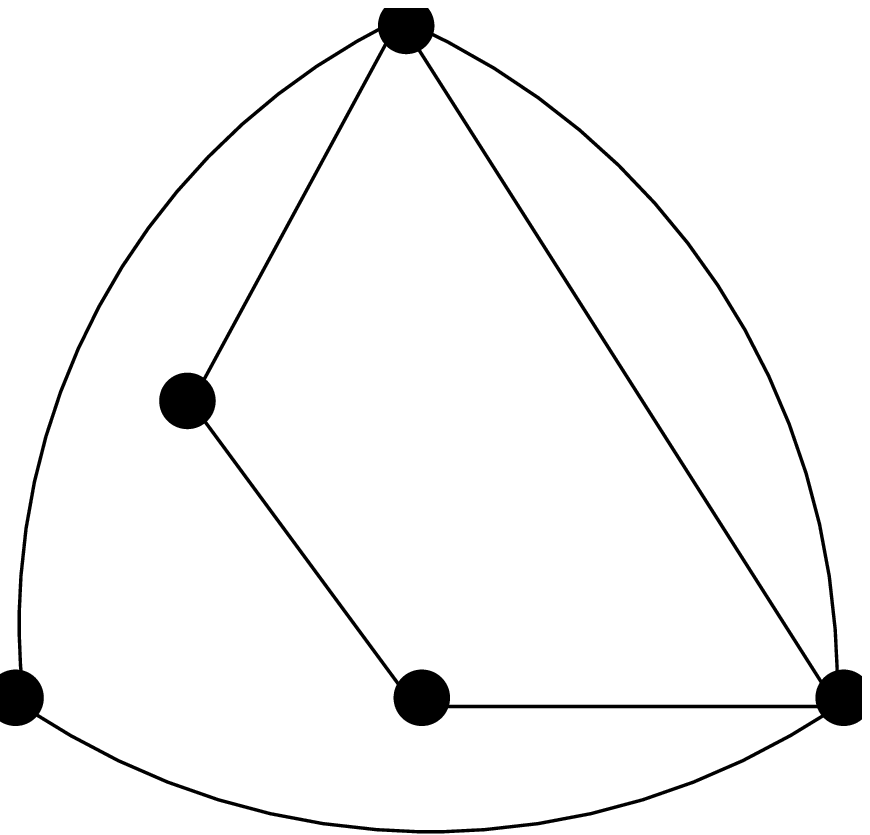}}

\end{center}
\end{proof}
\begin{proof}[Proof of the claim]
By contradiction, suppose that
$B_{\g}=\{0,1,2,b_1(\g)\}$. Since $b_1(\g)\in B_{\g}$, $\g$ is
eulerian by property (P\ref{h}). 
In particular, we can decompose $\g$ as an edge-disjoint
union of circuits: namely, 
 $\g=\gamma_1\cup\cdots\cup\gamma_r$ and 
$E(\gamma_i)\cap E(\gamma_j)=\emptyset$ for all $i,j$. 
Since there are only a finite number of possibilities
for such a decomposition of $\g$, 
we can choose a decomposition with $r$ maximal.

For all $i,j=1,\dotsc,r$, we have $\# V(\gamma_i)\cap
V(\gamma_j)\leq 1$ and $b_1(\gamma_i\cup\gamma_j)\leq 2$.
In fact, if $\gamma_i$ and $\gamma_j$ have $s\geq 3$ common vertices,
then $\gamma_i\cup\gamma_j$ has a decomposition as an edge-disjoint 
union of $s$ circuits, hence we get a decomposition of $\g$ in $r+s-2$ circuits,
against the maximality of $r$. Moreover, we have
\begin{align*}
 b_1(\gamma_i\cup\gamma_j)&=
\delta_{\gamma_i\cup\gamma_j}-\nu_{\gamma_i\cup\gamma_j}
+c(\gamma_i\cup\gamma_j)\\& = \delta_{\gamma_i}+\delta_{\gamma_j}-
\nu_{\gamma_i}-\nu_{\gamma_j}+\# V(\gamma_i)\cap V(\gamma_j)
+c(\gamma_i\cup\gamma_j)
\\&= \# V(\gamma_i)\cap V(\gamma_j) +c(\gamma_i\cup\gamma_j).
\end{align*} 
If $\# V(\gamma_i)\cap V(\gamma_j)=2$, then
$b_1(\gamma_i\cup\gamma_j)=3\in B_{\g}$,
a contradiction: so $\# V(\gamma_i)\cap
V(\gamma_j)\leq 1$ and $b_1(\gamma_i\cup\gamma_j)\leq 2$.

Consider now the cyclic set $\gamma_1\cup\cdots\cup\gamma_{r-1}$:
since $b_1(\gamma_1\cup\cdots\cup\gamma_{r-1})\in B_{\g}$, we get
\[ r-1\leq b_1(\gamma_1\cup\cdots\cup\gamma_{r-1})\leq 2, \]
thus $r\leq 3$. Since $r\leq 2$ would imply $b_1(\g)\leq 3$, we obtain
$r=3$. Then we have 
\begin{align*}
5&\leq b_1(\g)=\delta_{\g}-\nu_{\g}+1\\
&=\sum_{i=1}^3\delta_{\gamma_i}
-\sum_{i=1}^3\nu_{\gamma_i}+\sum_{1\leq i<j\leq 3}\#V(\gamma_i)\cap
 V(\gamma_j)-\#V(\gamma_1)\cap V(\gamma_2)\cap V(\gamma_3)+1\\
&\leq \sum_{1\leq i<j\leq 3}\#V(\gamma_i)\cap
 V(\gamma_j)+1\leq 4, \end{align*}
a contradiction.
\end{proof}

We conclude with a simple consequence involving stable spin curves 
\begin{cor}
Let $X$ be a stable curve of genus at least $4$, having superstable dual graph.
\begin{enumerate}[(i)]
\item If \  $2^{b_1(\g_X)-2}\not\in L(S_X)$, then $X$ is the union of two smooth, irreducible components meeting in
$b_1(\g_X)+1$ points.
\item If \  $2^{b_1(\g_X)-3}\not\in L(S_X)$,  and there exists $m\in L(S_X)$ such that $m<2^{b_1(\g_X)-3}$,
then $b_1(\g_X)=4$  and the dual graph of $X$ is a fat triangle.
\end{enumerate}
\end{cor}
Where the ``fat triangle" is defined in Theorem~\ref{tre}.
\begin{proof}  Just observe that $(i)$ and $(ii)$ translate, respectively,  Theorem~\ref{due} and 
  Theorem~\ref{tre}.
\end{proof}

\medskip

\small

\end{document}